\documentstyle{article}

\newcommand{\vs}{\vskip12mm}
\newcommand{\ms}{\vskip13pt}
\newcommand{\bs}{\vskip20pt}
\newcommand{\tem}[1]{\noindent {\bf Theorem {#1}.}}
\newcommand{\cor}[1]{\noindent {\bf Corollary {#1}.}}

\newcommand{\rrk}[1]{\noindent {\bf Remark {#1}.}}
\newcommand{\lem}[1]{\noindent {\bf Lemma {#1}.}}
\newcommand{\prop}[1]{\noindent {\bf Proposition {#1}.}}
\newcommand{\exam}[1]{\noindent {\bf Example {#1}.}}

\title{Finite groups with the same character tables, Drinfel'd algebras and
Galois algebras}

\author{A.A.Davydov}

\begin{document}

\maketitle

\begin{abstract}
We prove that finite groups have the same complex character tables
iff the
group algebras are twisted forms of each other as Drinfel'd
quasi-bialgebras or iff there is non-associative bi-Galois algebra over
these groups.
The interpretations of class-preserving automorphisms
and permutation representations with the same character
in terms of Drinfel'd algebras are also given.
\end{abstract}

\vs

\noindent {\bf 1.~Introduction.}
The theory of quasi-Hopf algebras was developed by V.G.Drinfel'd for the 
description of quantizations of Lie groups and algebras or so-called quantum
groups.

Althought the deformational quantization approach which is so useful
in the theory of quantum
groups can't be applied for the the case of finite groups,
the idea of twisting
seems to be very suitible for reformulating of various problems from 
representation theory of finite groups.

The key observations of this article is that any bijection between
character tables of
finite groups corresponds to the quasi-isomorphism of the group algebras 
considered as quasi-Hopf algebras and any two homomorphisms of the group 
algebras define the same map of character tables iff they are twisted forms.

In particular, we can give the definitions in terms of (quasi-)Hopf
algebras of
such objects as class-preserving automorphisms, permutation
representations with
the same character, groups with the same character tables.
Namely, any class-preserving automorphism is twisted form of identity maps as 
homomorphisms of Hopf algebras. Two permutation representations have the same 
complex character iff the corresponding homomorphisms into symmetric group are
twisted forms as homomorphisms of Hopf algebras. Two groups have the same 
character tables iff their group algebras are twisted forms as quasi-Hopf 
algebras. 

This point of view allows to select the subclass of pairs of groups with the 
same character tables. This subclass consists of pairs of groups whose group 
algebas are twisted forms as Hopf algebras. 

The notion of quasi-homomorphism of group algebras
can be formulated in terms of
Galois algebras. Using the calculation of automorphisms of associative Galois 
algebras we can describe quasi-isomorphisms of group algebras as Hopf algebras.
These quasi-isomorphisms correspond to normal abelean $2$-subgroups equipped with some
non-degenerated bimultiplicative forms.

\bs

\noindent {\bf 2.~Semirings and character tables.}
A {\em semiring} is a set $S$ with a collection of non-negative integers 
$\{ m^{x}_{x_1 ,x_2}, x,x_1 ,x_2\in S\}$ ({\em structural constants}) which 
satisfy the ({\em associativity}) condititon
$$m^{x}_{x_1 ,x_2 ,x_3} = \sum_{t\in S}m^{x}_{x_1 ,t}m^{t}_{x_2 ,x_3} = 
\sum_{s\in S}m^{s}_{x_1 ,x_2}m^{x}_{s,x_3}, \quad \forall x,x_1 ,x_2 ,
x_3\in S.$$
An element $e$ of the semiring $S$ is an {\em identity} if $m^{s}_{t,e} = 
m^{s}_{e,t} = \delta_{s,t}$ for all $s,t\in S$. 

A {\em morphism} from the semiring $S$ to the semiring $S'$ is a 
collection $\{ n^{s}_{t}, s\in S, t\in S' \}$ of non-negative integers which 
satisfy the following condition:
$$
\sum_{s\in S}m^{s}_{s_1 ,s_2}n^{t}_{s} = 
\sum_{t_1 ,t_2 \in S'}{m'}^{t}_{t_1 ,t_2}n^{t_1}_{s_1}n^{t_2}_{s_2} \eqno (1)
$$
for any $s_1 ,s_2 \in S$ and $t\in S'$. 
A {\em degree map} $d$ for the semiring $S$ is a morphism from $S$
to the one-element semiring with identity, e.g. a colection $\{ d(s), s\in S\}
$ of non-negative integers such that $d(s_1)d(s_2) =
\sum_{s\in S}m^{s}_{s_1 ,s_2}d(s)$.

The {\em enveloping ring} $A(S)$ of the semiring $S$ is the free 
$\bf Z$-module with the basis $\{[s], s\in S\}$ labeled by the elements of 
$S$ and with the multiplications $[i][j] = \sum_{s\in S}m^{s}_{i,j}[s]$. 
We will denote by $A_{\geq 0}(S)$ the cone of non-negative
combinations of basic elements (the cone of {\em non-negative elements}).

A morphism of semirings defines a homomorphism of their enveloping 
rings $f:A(S)\to A(S')$ where $f([s]) = \sum_{t\in S'}n^{t}_{s}[t]$. 

Denote by $(\ ,\ )$ canonical bilinear form on $A(S)$
$$(x,y) = \delta_{x,y},\quad\forall x,y\in S.$$
The semiring $S$ is {\em rigid} if there defined an antihomomorphism 
$(\ )^* :S\to S$ ({\em conjugation}) such that
$$(xy,z) = (y,x^* z),\quad\forall x,y,z\in S.$$
Note that conjugation is an anti-endomorphism of the enveloping ring $A(S)$:
$$(z,(xy)^*w) = ((xy)z,w) = (x(yz),w) = (yz,x^*w) = (z,y^*x^*w),$$
or $(xy)^* = y^*x^*$. 

It follows from the definition that the kernel of the conjugation $(\ )^*$
lies in the kernel of the bilinear form $(\ ,\ )$
$$(x,y) = (1,x^* y) = 0,\quad\mbox{for}\quad x\in ker(\ )^* , y\in S.$$
Since the bilinear form $(\ ,\ )$ is non-degenerated the conjugation is 
injective. So it is bijective for the finite semiring $S$. In that case
the conjugation has a finite order as an automorphism of the finite set $S$.

\ms

\lem{1}
{\it Let $d$ be a degree map for the rigid seniring $S$.
Then $\rho = \sum_{s\in S}d(s^*)s\in A(S)$ satisfies to the conditions
$$x\rho = d(x)\rho,\ \forall x\in A(S).$$
}

\noindent{\it Proof.}
Since $m_{x,s}^t = (t,xs) = (x^* t,s) = (t^* x,s^*) = m_{t^*,x}^{s^*}$ we have
$$x\rho = \sum_{s\in S}d(s^*)xs = \sum_{s,t\in S}d(s^*)m_{x,s}^t t = 
\sum_{s,t\in S}d(s^*)m_{t^*,x}^{s^*}t = \sum_{t\in S}d(t^* x)t = \rho\d(x).$$
$\Box$
\ms

{\bf Proposition 1}(Uniqueness of degree map){\bf.} 
{\it Any two degree maps for commutative rigid semiring $S$ coincides.
}

\noindent{\it Proof.}
Let $d,d'$ be degree maps for $S$. Define $\rho = \sum_{s\in S}d(s^*)s,
\rho' = \sum_{s\in S}d'(s^*)s$.
Then $d(\rho')\rho = \rho'\rho = d'(\rho)\rho'$, which means $d=d'$.
$\Box$
\ms

Since the enveloping ring $A(S)$ of rigid semiring $S$ is equipped with
non-degenerated semi-invariant bilinear form, the enveloping
algebra $A_{{\bf Q}}(S) = A(S)\otimes {\bf Q}$
over rational numbers ${\bf Q}$ is semisimple.

For commutative rigid semiring $S$ the enveloping algebra
$A_{{\bar {\bf Q}}}(S) = A(S)\otimes {{\bar {\bf Q}}}$ over
algebraic closure ${\bar {\bf Q}}$ of ${\bf Q}$ is isomorphic
to the algebra of functions on some finite set $Cl(S)$
({\em "conjugacy classes"} of $S$). The set $Cl(S)$ can be
identified with the set of maximal ideals of $A_{{\bar {\bf Q}}}(S)$,
so that the value $x(c)$ of $x\in A_{{\bar {\bf Q}}}(S)$ on $c\in Cl(S)$
is unique element of ${\bar {\bf Q}}$ such that $x\in x(c)+m_c$.
Here $m_c$ is maximal ideal of $A_{{\bar {\bf Q}}}(S)$ corresponding to $c$.

For any commutative rigid semiring $S$ we can associate a {\em character 
table} $(s(c))_{s\in S,c\in Cl(S)}$, which is $|S|\times|S|$-matrix
with entries in ${\bar {\bf Q}}$.

\ms

\prop{2}
{\it The map $f:A(S)\to A(S')$ given by the collection
$\{ n^{s}_{t}, s\in S, t\in S' \}$ is a homomorphism of (semi)rings
(the collection satisfies to the condition {\rm {(1)}}) if and only if
there is a map $f^*:Cl(S')\to Cl(S)$ such that $f(s)(c) = s(f^* (c))$
for any $s\in S, c\in Cl(S')$.
}

\noindent{\it Proof.}
The map $f:A(S)\to A(S')$ is a ring homomorphism iff
$f_{{\bar {\bf Q}}}:A_{{\bar {\bf Q}}}(S)\to A_{{\bar {\bf Q}}}(S')$
is a homomorphism of ${\bar {\bf Q}}$-algebras. Any homomorphism
of algebras of functions corresponds to the map of sets $f^*:Cl(S')\to Cl(S)$.
$\Box$
\ms

\exam{1} 
The set $Irr(G)$ of irreducible characters of the finite group $G$ has a 
natural semiring structure:
$$\chi\psi = \sum{\eta\in Irr(G)}m^{\eta}_{\chi ,\psi}\eta,\quad
\chi ,\psi\in
Irr(G).$$

The {\em map} of character tables of the finite groups $G_1, G_2$ is a pair
consisting of 

i) the map of the sets of conjugacy classes 
$cl(G_1 )\to cl(G_2),\quad C\mapsto C^*$,

ii) the map from the set of irreducible characters to the semiring of 
characters
$Irr(G_2 )\to R_{\geq 0}(G_1 ),\quad \chi\mapsto\chi^* =
sum_{\psi\in Irr(G_1 )}n_{\chi ,\psi}\psi$, where $n_{\chi ,\psi}\geq 0$

such that $\chi^* (C) = \chi (C^* )$ for all $\chi\in Irr(G_2 ),
C\in cl(G_1 )$.

We say that two finite groups $G_1 , G_2$ have {\em the same character tables} 
if there are one-to-one mappings $\chi\mapsto\chi^*$ and $C\mapsto C^*$
between
the sets of irreducible characters and conjugasy classes, respectively, of
$G_1$
and $G_2$, such that $\chi^* (C^* ) = \chi (C)$ for all $\chi , C$.
For examples, see [1, 3, 4, 8, 10, 12, 13].

\bs

\noindent {\bf 3.~Semisimple monoidal categories.}
The {\em monoidal category} [6] is a category
{\cal G} \ with a bifunctor
\begin{displaymath}
\otimes :{\cal G} \times {\cal G} \longrightarrow
{\cal G} \qquad (X,Y) \mapsto
X \otimes Y
\end{displaymath}
which called {\em tensor (or monoidal) product}. This functor must be equiped 
with a functorial collection of isomorphisms (so-called {\em associativity 
constraint})
\begin{displaymath}
\varphi_{X,Y,Z} : X \otimes (Y \otimes Z) \rightarrow (X \otimes Y) \otimes Z 
\qquad \mbox{for any} \quad X,Y,Z \in {\cal G}
\end{displaymath}
which satisfies to the following {\em pentagon axiom}:
$$(X\otimes\varphi_{Y,Z,W})\varphi_{X,Y\otimes Z,W}(\varphi_{X,Y,Z}\otimes W) = 
\varphi_{X,Y,Z\otimes W}\varphi_{X\otimes Y,Z,W}.$$

A {\em quasimonoidal functor} between monoidal categories ${\cal G}$
and ${\cal H}$ is a functor $F : {\cal G} \longrightarrow
{\cal H}$ , which is equipped with the functorial collection
of isomorphisms (the so-called {\em quasimonoidal structure})
\begin{displaymath}
F_{X,Y} : F(X \otimes Y) \rightarrow F(X) \otimes F(Y)
\qquad \mbox{for any} \quad X,Y \in {\cal G}. 
\end{displaymath}
We shal call it {\em monoidal structure} if 
$$F_{X,Y\otimes Z}(I\otimes F_{Y,Z}) = F_{X\otimes Y,Z}(F_{X,Y}\otimes I)$$
for any objects $X,Y,Z \in {\cal G}$.

The morphism $\alpha:F\to G$ between two monoidal functors
$F,G:{\cal G}\to{\cal H}$ is {\em monoidal} if
$F_{X,Y}(\alpha_X\otimes\alpha_Y) = \alpha_{x\otimes Y}G_{X,Y}$ for any
$X,Y\in{\cal G}$.

Monoidal category structures on ${\cal G}$ differed by the associativity 
constraint will be called {\em twisted forms} of each other. 

The structures of monoidal functor for $F: {\cal G}\to {\cal H}$ 
will be called {\em twisted forms} of each other.  

The monoidal category ${\cal G}$ is {\em rigid} if it is equipped
with the {\em dualization} functor, which is a contravariant functor
$(\ )^* :{\cal G}\to{\cal G}$ with a collections of morphisms
$\iota:1\to X\otimes X^*$ and $\epsilon\upsilon:x^*\otimes X\to 1$ for any
$X\in{\cal G}$ such that the compositions
$$X\stackrel{I\otimes\iota}{\longrightarrow}X\otimes (X^*\otimes X)
\stackrel{\varphi}{\longrightarrow}(X\otimes X^*)\otimes X
\stackrel{\epsilon\upsilon\otimes I}{\longrightarrow}X,$$
$$X^*\stackrel{\iota\otimes I}{\longrightarrow}(X^*\otimes X)\otimes X^*
\stackrel{\varphi^{-1}}{\longrightarrow}X^*\otimes (X\otimes X^*)
\stackrel{I\otimes\epsilon\upsilon}{\longrightarrow}X^*$$
are identical. 

Let ${\cal G}$ be a semisimple monoidal $k$-linear category over
the field algebraically closed $k$
with the set $S$ of isomorphism classes of simple objects. The collection of
dimensions $m^{x}_{y,z} = dim_k Hom_{\cal G}(X,Y\otimes Z)$ form a semiring 
structure on the set $S$. Here $X,Y$ and $Z$ are some representatives of the 
classes $x,y,z\in S$. Note that the enveloping ring of semiring $S$
coincides with the Grothendieck ring $K_0 ({\cal G})$ of the category
${\cal G}$.
The semiring $S({\cal G})$ is rigid for the rigid monoidal category ${\cal G}$.

\ms

\prop{3}
{\it Semisimple monoidal categories are twisted forms of each other iff
their semirings of simple objects are isomorphic.
Isomorphism classes of quasimonoidal functors
$F:{\cal G}\to{\cal H}$ between semisimple monoidal categories
are in one-to-one correspondence with the homomorphisms
$S({\cal G})\to S({\cal H})$ of the semirings of simple objects.
In particular, monoidal functors $F,G:{\cal G}\to{\cal H}$
between semisimple monoidal categories are twisted forms
of each other iff they induce the same map
$K_0 ({\cal G})\to K_0 ({\cal H})$ of the Grothendieck rings.
}

\noindent{\it Proof.}
The proposition follows from the fact that two fuctors
$F,G:{\cal G}\to{\cal H}$ between semisimple categories are isomorphic iff
they induce the same map $S({\cal G})\to S({\cal H})$ of the semirings
of simple objects.
$\Box$
\bs

\noindent {\bf 4.~Drinfel'd algebras.}
A {\em Drinfel'd algebra} or {\em quasi-bialgebra} [7]
is an algebra $H$ together
with a homomoprhisms of algebras
$$\Delta :H\rightarrow H\otimes H\quad\mbox{(coproduct)},
\quad\varepsilon :H\rightarrow k\quad\mbox{(counit)}$$
and an invertible element $\Phi\in H^{\otimes 3}$ ({\em associator})
for which
$$(\Delta\otimes I)(\Delta (h)) = \Phi (I\otimes\Delta)(\Delta (h))\Phi^{-1}
\quad\forall h\in H\quad\mbox{(coassociativity)},$$
$$(I\otimes I\otimes\Delta)(\Phi)(\Delta\otimes I\otimes I)(\Phi) = 
(1\otimes\Phi)(I\otimes\Delta\otimes I)(\Phi)(\Phi\otimes 1),$$
$$(\varepsilon\otimes I)\Delta = (I\otimes\varepsilon)\Delta = I.$$

Drinfel'd algebra is a generalization of the well-known notion of 
{\em bialgebra} which corresponds to the case of trivial associator
$\Phi = 1$.

Drinfel'd algebras structures on the algebra $H$
which is differed only by associator will be called
{\em twisted forms} of each other.

A {\em quasi-homomorphism} of quasi-bialgebras $H_1$ and $H_2$ is pair
$(f,F)$ consisting of a homomorphism of algebras $f:H_1\to H_2$
and an invertible element $F\in H_{2}^{\otimes 2}$ such that
$$\Delta (f(h)) = F(f\otimes f)(\Delta (h))F^{-1}.$$
It is a {\em homomorphism} of quasi-bialgebras if, additionly, 
$$(\Delta\otimes I)(F)(F\otimes 1)(f\otimes f\otimes f)(\Phi_1 ) = 
\Phi_2 (I\otimes\Delta )(F)(1\otimes F).$$

Two homomorphisms of quasi-bialgebras are {\em twisted forms} of each other
if they differ only by the invertible element $F$.
We can define the {\em morphism} between two homomorphisms
$(f,H),(g,G):H_1\to H_2$ as an element $c\in H_2$ such that $cf(h) = g(h)c$
for any $h\in H_1$ and $\Delta(c)G = F(c\otimes c)$.

A {\em homomorphism of bialgebras} $H_1 ,H_2$ is a 
homomorphism of algebras $f:H_1\to H_2$ such that
$\Delta f = (f\otimes f)\Delta .$

Now we discuss the connection between monoidal categories
and quasi-bialgebras.
Coproduct allows to define the structure of $H$-module on the tensor product 
$M\otimes_{k}N$ of two $H$-modules: 
$$h*(m\otimes n) = \Delta(h)(m\otimes n),\qquad h\in H, m\in M, n\in N.$$ 
The associator $\Phi$ defines the associativity constraint 
$$\varphi :L\otimes M\otimes N\rightarrow L\otimes M\otimes N, \qquad 
\varphi (l\otimes m\otimes n) = \Phi(l\otimes m\otimes n).$$
Thus the category ${\cal M} (H)$ of $H$-modules is a monoidal category. 
\newline
The homomorphism of quasi-bialgebras $f:H_1\to H_2$ defines the monoidal 
functor 
$$f^* :{\cal M} (H_{2})\to{\cal M} (H_{1})$$ 
with the monoidal structure defined by the invertible element 
$F\in H_{2}^{\otimes 2}$ 
$${f^*}_{M,N}:f^* (M\otimes N)\to f^* (M)\otimes f^* (M)\quad
{f^*}_{M,N}(m\otimes n) = F(m\otimes n).$$
The morphisms between homomorphisms $f,g:H_1\to H_2$
of quasi-bialgebras correspond to the monoidal morphisms
between monoidal functors $f^*,g^* :{\cal M} (H_{2})\to{\cal M} (H_{1})$.

The quasi-Hopf algebra $H$ will be called {\em rigid}
if the monoidal category ${\cal M} (H)$ is rigid.
The dualization functor for ${\cal M} (H)$ corresponds to the
antihomomorphism $S:H\to H$ ({\em antipode}) with
some additional properties (see [7]).
For bialgebra these properties takes a form of the relation
$$\mu(S\otimes I)\Delta = \mu(I\otimes S)\Delta = \varepsilon,$$
where $\mu :H\otimes H\to H$ is the multiplication in $H$.
Bialgebra with an antipode is called {\em Hopf algebra}.

\ms

\exam{2}
Group algebra $k[G]$ of the group $G$.
As $k$-vector space $k[G]$ is spanned by the elements of the group $G$.
The structure maps have the following forms on the basis: 
$$\Delta (g) = g\otimes g,\quad \varepsilon (g) = 1,\quad S(g) = g^{-1}.$$

\ms

The homomorphism of the groups $f:G_1\to G_2$ defines
the homomorphism of their group algebras and any homomorphism
of bialgebras $k[G_1 ]\to k[G_2 ]$ is of this kind.
The group algebra provides an example of so-called
{\em cocommutative} Hopf algebra $t\Delta = \Delta$.
Over the algebraically closed field of characteristic zero
group algebras are characterized by this property
(Kostant theorem): any cocommutative finite dimensional Hopf algebra
is isomorphic to the group algebra.

For semisimple quasi-bialgebra $H$ denote by $S(H)=S({\cal M} (H))$
the semiring of simple modules.
The semiring $S(H)$ is rigid for quasi-Hopf algebra $H$. 

The next proposition is the direct consequence of the
definitions and proposition Proposition~ 3.

\ms

\prop{4}
{\it The homomorphisms of quasi-bialgebras $f_1 ,f_2 :H_1\to H_2$
are twisted forms if and only if the monoidal functors $(f_1 )^* ,(f_2 )^*$
are twisted forms.
In particular, the homomorphisms of semisimple quasi-bialgebras
$f_1 ,f_2 :H_1\to H_2$ induce the same homomorphism
$K_0 (f_1) ,K_0 (f_2):K_0 (H_2)\to K_0 (H_1)$
of Grothendieck rings if and only if one is isomorphic to the twisted form
of the other.
}

\ms

The generalization of the so-called Tannaka-Krein theory [5, 6]
states that we can reconstruct a quasi-bialgebra
from the monoidal category ${\cal G}$ and a quasimonoidal functor
$F:{\cal G}\to{\cal M} (k)$ to the category of vector spaces
as endomorphisms $End(F)$ of the functor $F$.

\ms

\tem{1}
{\it Finite dimensional semisimple quasi-Hopf algebras $H_1 ,H_2$
are quasi-isomorphic if and only if their semirings of simple objects
$S(H_2),S(H_1)$ are isomorphic.
}

\noindent{\it Proof.}
Since twisting does not change the semiring of simple objects
we need to prove the if statement. Let $f^*:S(H_2)\to S(H_1)$
be an isomorphism of semirimgs of simple objects. By Proposition~ 3 we can costruct a quasi-monoidal functor
(equivalence) $F:{\cal M} (H_2)\to{\cal M} (H_1)$
which induces the given homomorphisms $f^*$.
By Proposition~ 1 the composition $d_1 f^*$
coincides with $d_2$, where $d_i$ is a degree map for $S(H_i)$.
Hence the composition $F_1 F$ of functor $F$ with the forgetful funtor
$F_1:{\cal M} (H_1)\to{\cal M} (k)$ is isomorphic to the forgetful funtor
$F_2:{\cal M} (H_2)\to{\cal M} (k)$ as quasi-monoidal functor. Using Tannaka-Krein theory we can construct
the isomomorphism of quasi-bialgebras $f:H_1\to H_2$ as
$$H_1 = End(F_1)\to End(F_1 F)\to End(F_2) = H_2.$$
$\Box$
\ms

\rrk{1}
The weak version of the theorem Theorem~ 1 for Hopf algebras
was proved in [11] where it was assumed that the isomorphism of
(enveloping algebras of) semiring preserves the class
of regular representation.

\ms

\cor{1}
{\it The finite groups $G_1 ,G_2$ have the same character table
if and only if their group algebras are isomorphic as quasi-Hopf algebras,
e.g. there is an isomorphism of algebras $f:k[G_1 ]\to k[G_2 ]$
and an invertible element $F\in k[G_2 ]^{\otimes 2}$ such that
$F\Delta_2 (f(x)) = (f\otimes f)(\Delta_1 (x))$ for any $x\in k[G_1 ]$.
}

\ms

If we denote by $\Delta_F$ the twisted by $F$ comultiplication on $k[G_2 ]$ 
$\Delta_F (x) = F\Delta_2 (x)F^{-1}$ then the map $f$ would be
an isomorphism of Hopf algebras
$(k[G_1 ],\Delta_1 )$ and $(k[G_2 ],\Delta_F )$.
The existence of such isomorphism is equivalent to the existence
of an isomorphism of groups
$$G_1\to G(F) = G(k[G_2 ],\Delta_F ) = \{ x\in k[G_2 ], F\Delta_2 (x) =
(x\otimes x)F\}.$$
The cocommutativity of the coproduct $\Delta_F$
is equivalent to the condition
$$
t(F) = F \eqno (2)
$$
The coassociativity of the twisted coproduct $\Delta_F$
is equivalent to the equation on $F$
$$
(1\otimes F)(I\otimes\Delta )(F) = (F\otimes 1)(\Delta\otimes I)(F)\Phi,
\eqno (3)
$$
where $\Phi$ is some invertible $G_2$-invariant element of
$k[G_2 ]^{\otimes 3}$ ({\em associator}).
In particular, such $\Phi$ satisfy to the equation
$$
(\Phi\otimes 1)(I\otimes\Delta\otimes I)(\Phi )(1\otimes\Phi ) =
(\Delta\otimes I\otimes I)(\Phi )(I\otimes I\otimes\Delta )(\Phi ). \eqno (4)
$$
We can replace $F$ by the product $FC$ for any $G_2$-invariant element
$C\in k[G_2 ]^{\otimes 2}$ without changing the twisted coproduct
$\Delta_{FC} = \Delta_F$. The $G_2$-invariance of $C$ allows
to write the associator $\Phi^C$ for the product $FC$ as
$$
\Phi^C = (\Delta\otimes I)(C)^{-1}(C\otimes 1)^{-1}\Phi
(1\otimes C)(I\otimes\Delta )(C). \eqno (5)
$$
Thus the element $\Phi$ is defined up to the transformations (5). 

We can also replace $F$ by $F^c = (c\otimes c)F\Delta(c)^{-1}$
for invertible $c\in k[G_2 ]$.
Then the corresponding twisted coproducts will be connected
by conjugation by $c$
$$\Delta_{F^c}(cxc^{-1}) =  (c\otimes c)\Delta_F (x)(c\otimes c)^{-1}.$$

The preveous theorem reduces the problem of finding
finite groups whose character tables coincide with the character table of
$G$ to the problem of
finding the solutions $(\Phi ,F)$
to the equations (2), (3), (4)
such that the order of the group $G(F)$ equals $|G|$.
If the ground field $k$ is algebraically closed of characteristics zero,
then we can ommite the condition $|G(F)| = |G|$ using Kostant theorem.

In [7] V.Drinfeld suggested the following way of solving
the equation (3) for general Hopf algebra. Introduce the new multiplication
$\mu_F$ on the dual Hopf algebra $k(G) = k[G]^*$
$$\mu_F (l\otimes m)(x) = (l\otimes m)(F\Delta (x)),\quad
l,m\in k(G), x\in k[G].$$
This multiplication will be invariant under the action of the group
$G$ on $k(G)$
$$(gl)(x) = l(xg).$$
By another words, elements of the group $G$ act as automorphisms
of the algebra $R_F = (k(G),\mu_F )$. Moreover, the algebra
$R_F$ is a so-called {\em Galois} $G$-algebra.
It means, that the natural map of vector spaces
$$R_F \otimes R_F \to Hom(k[G],R_F ),\quad
l\otimes m\mapsto (g\mapsto g(l)m )$$
is an isomorphism. The group $G(F)$
appears as automorphisms group $Aut_G (R_F )$ of $G$-algebra $R_F$.
It is not hard to see that if $|G(F)| = |G|$, then $R_F$ is also Galois
$G(F)$-algebra, or bi-Galois $G-G(F)$-algebra.

The algebras $R_{F_1},R_{F_2}$ are isomorphic as $G$-algebras iff
there is an invertible $c\in k[G]$ such that $F_1 = F_{2}^{c}$.
This method is mostly applicable for the case of $\Phi = 1$,
because of the following fact:

the algebra $R_F$ is associative iff $\Phi = 1$. 

In general, it would be only {\em $\Phi$-associative} in the following sense:
$$x(yz) = \Phi (xy)z,\quad \forall x,y,z\in R,$$
where the product
$\Phi (xy)z = \mu (\mu\otimes I)(\Phi (x\otimes y\otimes z))$
is defined by the action of $k[G]^{\otimes 3}$ on $R^{\otimes 3}$.

\bs

\noindent {\bf 5.~Galois algebras.}
Here we give brief description of bi-Galois associative algebras. As was explained above they correspond to the isomorphisms of character tables with trivial associators. 

Let $R$ be an algebra with the action of the group $G$ ({\em $G$-algebra}).
The {\em cross product} $R*G$ is a vector space spanned by the elements 
$a*g\quad a\in R, g\in G$ satisfying $(a+b)*g = a*g + b*g$.
The multiplication is given by the formula
$$(a*g)(b*f) = ag(b)*gf\qquad\forall a,b\in R, g\in G.$$

A $G$-algebra $R$ is {\em Galois} if the map 
$$\theta: R*G\to End(R)\qquad\theta (a*g)(b) = ag(b)$$
is an isomorphism. 

A Galois $G$-algebra $R$ has the following properties:

$R$ has no non-trivial $G$-invariant twosided ideals,

$R$ is semisimple,

$G$ acts transitively on the set of maximal twosided ideals of $R$. 

Let $S$ be a stabilizer of some maximal ideal $M$ of $R$. Then the quotient 
algebra $B=R/M$ is simple Galois $S$-algebra and $R$ can be identified with 
the algebra of $S$-invariant functions
$$ind_{S}^{G}(B) = 
\{a:G\to B,\quad a(sg) = s(a(g))\qquad\forall s\in S, g\in G\}$$
with the $G$-action $(fa)(g) = a(gf)$.

The $S$-algebra $B = End(V)$ is Galois iff the multiplier of the projective
representation $S\to PGL(V) = Aut(B)$ is a non-degenerated 2-cocycle.
We call a 2-cocycle $\alpha\in Z^2 (G,k^* )$ {\em non-degenerated} if
for any $s\in S$ the map from the centralizer 
$$C_S (s)\to k^*\qquad t\mapsto \alpha (s,t)\alpha (t,s)^{-1}$$
is non-trivial. 

\ms

\exam{3}
Let $A$ be a finite abelian group. Denote by $\hat A$ the dual group 
$Hom(A,k^*)$. The 2-cocycle $\alpha$ on $S = A\oplus\hat A$
$$\alpha ((a,\chi),(b,\psi)) = \chi(b),\quad a,b\in A,\chi ,\psi\in\hat A$$
is non-degenerated.

\ms

Describe the automorphisms of Galois algebras. 

The set of maximal ideals of $G$-Galois algebra $R$ can be identified as
$G$-set with $G/S$ where $S$ is a stabilizer of some maximal ideal. 

The action of automorphisms on maximal ideals defines the homomorphism
$Aut_G (R)\to N_G (S)/S$. 

The kernel of this homomorphism coincides with $Aut_S (B)$. 
The group of automorphisms of the simple Galois $S$-algebra is
isomorphic to the character group $\hat S = Hom_{group}(S,k^* )$.

The image of this homomorphism coincides with the stabilizer 
$St_{N_G (S)/S}(\alpha )$ of the cohomological class $\alpha\in H^2 (S,k^* )$. 

Thus we have a short exact sequence
$$\hat S\to Aut_G (R)\to St_{N_G (S)/S}(\alpha ).$$

The class of this extension in $H^2 (St_{N_G (S)/S}(\alpha ),\hat S)$ 
is the image of the class $\alpha\in H^2 (S,k^* )$ by 
$$d_{2}^{0,2}:H^0 (St_{N_G (S)/S}(\alpha ),H^2 (S,k^* ))\to  
H^2 (St_{N_G (S)/S}(\alpha ),H^1 (S,k^*))$$ 
the differential of 
Hochschild-Serre spectral sequence corresponding to the extension 
$S\to N_G (S)\to N_G (S)/S$. 

We can apply now the outlined description of automorphisms of
Galois algebras to the investigation of bi-Galois algebras.

A {\em biGalois} $G_1 -G_2$-algebra is an algebra $R$ with the commuting actions 
of the groups $G_1 ,G_2$ such that $R$ is both Galois $G_1$-algebra and 
$G_2$-algebra. 

The $G_1 -G_2$-biGalois algebra corresponds to the following data:

the normal inclusions $S\to G_i$ of the abelian group $S$
with the same quotient group $F$,

the non-degenerated class $\alpha\in H^2 (S,k^* )$ such that 
$$d(\alpha) = \gamma_1 - \gamma_2 ,$$
where $\gamma_i$ is the class of the extension $S\to G_i\to F$ in $H^2 (F,S)$
and $d$ is the differential of Hochschild-Serre spectral sequence
corresponding to the splitting extension of $F$ by $S$.

Functoriality of the differential $d$ allows to reduce consideration
to the case of $p$-group
$S$. The diffirential $d$ is trivial for abelian $p$-groups if $p\not= 2$.

\ms

{\bf Example} (see, also $\rm [8]$){\bf.}
Let $S$ be $2n$-dimensional vector space over the field ${\bf F}_2$ of two 
elements. The standard symplectic form $<\ ,\ >$ on $S$ defines a 
non-degenerated 2-cocycle 
$$\alpha\in Z^2 (S,k^* ),\quad\alpha (s,t) = (-1)^{\beta (s,t)},$$
where $\beta$ is bilinear form on $S$ such that
$<s,t>=\beta (s,t)-\beta (t,s)$.
\newline
Let $F=Sp_{2n}(2)$ be the group of automorphisms of the form $<\ ,\ >$.
For $n>1$ the cohomology group
$H^2 (F,\hat S) = H^2 (Sp_{2n},{\bf F}_{2}^{2n})$
is one dimensional ${\bf F}_2$-vector space generated by the class
$d(\alpha )$.
Thus the affine symplectic group $AffSp_{2n}(2)$ and the (unique) non-split 
extension of $Sp_{2n}(2)$ by ${\bf F}_{2}^{2n}$ have the same
character tables.
\newline
For $n=1$ the cohomology group $H^2 (Sp_{2n},{\bf F}_{2}^{2n})$
is trivial and the pair $(S,\alpha )$ defines the automorphism
of character table of
$AffSp_{2}(2)=S_4$ which doesn't correspond to any group automorphism.
This isomorphism intertwings the characters $\chi_4$ and $\chi_5$
and the classes $2A$ and $4A$.

\ms

\begin{tabular}{r|rrrrr} 
$S_4$    & 1 & 2A & 2B & 3A & 4A \\
\hline
$\chi_1$ & 1 &  1 &  1 &  1 &  1 \\
$\chi_2$ & 1 & -1 &  1 &  1 & -1 \\
$\chi_3$ & 2 &  0 &  2 & -1 &  0 \\
$\chi_4$ & 3 &  1 & -1 &  0 & -1 \\
$\chi_5$ & 3 & -1 & -1 &  0 &  1 \\
\end{tabular}

\bs

\noindent {\bf 6.~Class-preserving automorphisms and permutation
representations with the same character.}
We will call an automorphism $\phi\in Aut(G)$ of the finite group $G$
by {\em class-preserving} if $\phi$ preserves all conjugacy
classes of $G$ $\phi (g)\in g^G$ for all $g\in G$ (see [2, 9, 14]).

\ms

\prop{5}
{\it For any class-preserving automorphism $\phi$ of the finite group $G$
there is an invertible element $c\in k[G]$ such that
$\phi (g) = cgc^{-1}$ for any $g\in G$ and
$F = \Delta (c)^{-1}(c\otimes c)$ is $G$-invariant element of
$k[G]^{\otimes 2}$.
}

\noindent{\it Proof.}
It follows directly from the Proposition~ 4 and the fact that
class-preserving automorphism induces trivial automorphism
of the character ring $R(G)$.
$\Box$
\ms

{\em Permutation representation} of the group $G$ is a homomorphism
$\phi :G\to S_n$ to the group of automorphisms $S_n = Aut(X)$
of the finite set of order $|X|=n$.
Define the {\em character} $\chi_\phi$ of the permutation representation
$\phi :G\to S_n$  as $\chi_{\phi}(g) = |\{ x\in X, \phi (g)(x)=x\}|$.
So $\chi_\phi$ is an image of the natural $n$-dimensional
character of $S_n$ under the homomorphism $\phi^* :R(S_n )\to R(G)$.

\ms

\prop{6}
{\it For any two permutation representations $\phi ,\psi :G\to S_n$
of the finite group $G$ there is an invertible element $c\in k[S_n ]$
such that $\phi (g) = c\psi (g)c^{-1}$ for any $g\in G$ and
$F = \Delta (c)^{-1}(c\otimes c)$ is $\psi (G)$-invariant element of
$k[G]^{\otimes 2}$.
}

\noindent{\it Proof.}
It can be proved that the homomorphisms $\phi^* ,\psi^* :R(S_n )\to R(G)$
of character tables corresponded to permutation representations
$\phi ,\psi :G\to S_n$ coincides if coincides the characters
$\chi_\phi$, $\chi_\psi$.
Hence we can apply the Proposition~ 4. 
$\Box$
\bs

\noindent {\bf 7.~Concluding remarks.}
The cohomological nature of the sets of possible twistings
was actively explored
in the theory of quantum groups. Nonabelianity of those cohomology
is probabily
a major difficulty of the theory. In quantum group theory this difficulty
was overcome by methods of tangent cohomology which are unapplicable
for finite
groups. In this case non-abelian cohomology sets of twistings
can be abelianized
by means of algebraic K-theory. Namely, the maps from
the sets of twistings to
some Hochschild cohomology of representation ring can be costructed
[5]. The detailed description of those maps would be the subject
of subsequent paper.

\end{document}